\documentclass[11pt,english]{amsart}
\usepackage[T1]{fontenc}
\usepackage[latin9]{inputenc}
\usepackage{amstext}
\usepackage{amsthm}

\makeatletter
\usepackage{amsmath,amsfonts,amssymb,amsthm,epsfig}

\usepackage{color}
\usepackage[final,allcolors=blue,colorlinks=true]{hyperref}

\def\al{\alpha}
\def\be{\beta}
\def\ga{\gamma}

\def\ta{\theta}

\def\na{\nabla}
\def\Ga{\Gamma}  
\def\Om{\Omega}  

\def\cal{\mathcal}
\def\wq{\infty}

\def\loc{\text{\rm loc}}

\def\mf{\mathfrak}
\def\g{\mathfrak{g}}

\newcommand{\medint}{-\kern -,375cm\int}         
\newcommand{\medintinrigo}{-\kern -,315cm\int}
\newcommand{\wto}{\rightharpoonup}                

\newcommand{\YM}{\operatorname{YM}}

\numberwithin{equation}{section}
\newtheorem{theorem}{Theorem}[section]
\newtheorem*{theorem*}{Theorem}  

\newtheorem*{conclusion*}{Conclusin}

\newtheorem*{corollary*}{Corollary}

\newtheorem{definition}[theorem]{Definition}

\newtheorem{lemma}[theorem]{Lemma}
\newtheorem*{lemma*}{Lemma}

\newtheorem*{notation*}{Notation}

\newtheorem*{proposition*}{Proposition}

\newtheorem*{remark*}{Remark}

\newtheorem*{example*}{Example}                
\theoremstyle{definition}

\makeatother

\usepackage{babel}
\begin{document}
\title[Weak compactness of $\Omega$-Yang-Mills connections]{A note on weak compactness of $\Omega$-Yang-Mills connections} 

    \author[C.-Y. Guo and C.-L. Xiang]{Chang-Yu Guo and Chang-Lin Xiang*}         

\address[Chang-Yu Guo]{Research Center for Mathematics, Shandong University, Qingdao, P. R. China and Frontiers Science Center for Nonlinear Expectations, Ministry of Education, P. R. China and Department of Physics and Mathematics, University of Eastern Finland, 80101, Joensuu, Finland}
\email{guocybnu@gmail.com}

\address[Chang-Lin Xiang]{Three Gorges Mathematical Research Center, China Three Gorges University,  443002, Yichang,  P. R. China}
\email{changlin.xiang@ctgu.edu.cn}

\thanks{C.-Y. Guo is supported by  the Young Scientist Program of the Ministry of Science and Technology of China (No.~2021YFA1002200), the National Natural Science Foundation of China (No.~12311530037), the Taishan Scholar Project and the Natural Science Foundation of Shandong Province (No.~ZR2022YQ01). The corresponding author C.-L. Xiang was financially supported by NSFC (No. 12271296) and NSF of Hubei Province (No. 2024AFA061).}

\begin{abstract}
In this note, applying a compensation compactness argument developped by Chen and Giron (arXiv.2108.13529) on Yang-Mills fields, we extends their weak continuity result to the more general class of $\Omega$-Yang-Mills connections on principle bundles over compact Riemannian manifold.
\end{abstract}

\maketitle

{\small    
\smallskip
\newline
}
\bigskip

\section{Introduction and main result}

Let $M_{g}^{n}$ be an oriented closed (compact without boundary)
Riemannian manifold of dimension $n\ge4$ and $P\to M$ a principle
$G$-bundle with compact Lie group $G$. Let $\Om$ be a smooth $n-4$
form on $M$. Throughout we assume that $G$ is a matrix group
and the corresponding Lie algebra $\mf{g}$ is equipped with a $G$-equivariant
inner product $\langle,\rangle_{\mf{g}}$. In this note we aim to
study weak compactness for sequences of uniformly bounded weak $\Om$-Yang-Mills connections
on $P$.

$\Om$-Yang-Mills connections (or $\Om$-YM connections, for short)
are defined as critical points of the $\Om$-Yang-Mills functional
\[
\YM_{\Om}(A)=\int_{M}|F_{A}|^{2}dV-\int_{M}{\rm tr}(F_{A}\wedge F_{A})\wedge\Om
\]
over all connection 1-forms $A$ on $P$, where $F_{A}$ is the curvature
form of $A$. The Euler-Lagrange equation of which is given by
\[
D_{A}^{\ast}\left(F_{A}+\ast(F_{A}\wedge\Om)\right)=0.
\]
In this case, we also call $F_{A}$ an $\Om$-Yang-Mills field. An
important example of an $\Om$-YM field is the so-called $\Omega$-anti-selfdual
Yang-Mills instantons, which are defined as solutions of the equation
\[
\ast F_{A}=-\Om\wedge F_{A}.
\]
Such connection is not necessarily Yang-Mills unless $\Om$ is closed.
Another important example is the Hermitian-Yang-Mills connections
over general complex manifolds. Clearly, a $0$-Yang-Mills connection
is the standard Yang-Mills connection. In fact, if $\Om$
is a closed $n-4$ form, then the quantity $\int_{M}{\rm tr}(F_{A}\wedge F_{A})\wedge\Om$
is topological for $M$, and so critical points of $\YM_{\Om}$ are
identical to standard Yang-Mills connections. For more background
of $\Om$-YM theory, we refer to e.g. \cite{CDFN-1983-NP,chen-Wentworth-2022-CV,Donaldson-Segal-2011,Donaldson-Thomas-1998,Tian-2000-Annals}. 

In this note we are concerned with the problem of weak compactness
of $\Om$-Yang-Mills connections. That is, to show that the  limit of a sequence of weakly convergent $\Om$-Yang-Mills connections is a $\Om$-Yang-Mills connection as well.   Such kind of  compactness problem
is often of basic importance in the study of many important geometric
variational problems, such as Yang-Mills equations, harmonic maps,
wave maps and isometric immersions, see e.g. \cite{Chen-Li-2018-JGA,Chen-Li-2021-ARMA,chen-Wentworth-2022-CV,Donaldson-Kronheimer-book,Freire-M-S-1997,Freire-M-S-1998,Nakajima-1988,Uhlenbeck-1982-CMP}
and the references therein.   In a recent work \cite{chen-Wentworth-2022-CV},
Chen and Wentworth proved that if $\{A_{i}\}$ is a sequence of \emph{smooth
$\Om$-YM connections} with $\|F_{A_{i}}\|_{L^{2}}$ uniformly bounded,
then away from an $(n-4)$-rectifiable set in $M$, $A_{i}$ converges
locally uniformly in $C^{\wq}$-topology to a smooth $\Om$-YM connection
$A_{\wq}$ modulo gauge transformations; see also \cite[Theorem 5]{chen-Wentworth-2022-CV}
for more compactness results on smooth $\Om$-YM connections.  
 
In an interesting recent work \cite{Chen-Giron-21}, Chen and Giron proved
that for any dimension $n\ge4$, if $A_{i}\in L_{\loc}^{4}(P)$ is
a sequence of \emph{weak Yang-Mills connections} which converges weakly in
$L_{\loc}^{4}$ to a weak connection $A$, and $\|F_{A_{i}}\|_{L^{2}}$
is locally uniformly bounded, then $A$ is a also weak Yang-Mills
connection. This is often called \emph{weak compactness} or \emph{weak continuity}
of Yang-Mills connections in literature. Motivated by this interesting work, we
aim to establish a similar weak compactness result for sequences of
\emph{weak $\Om$-YM connections}. In fact, Chen and Wentworth \cite{chen-Wentworth-2022-CV} showed that $\Om$-YM connections have monotonicity properties and the singularity of $\Om$-YM connections can be  stratified according to symmetry. In view of the recent groudbreaking work of Naber and Valtorta \cite{Naber-Val-2017-Ann, Naber-Val-2019-Invent} on the singularity theory of harmonic mappings and energy identity on Yang-Mills connections, it is natural to expect that similar theories for $\Om$-YM connections hold. To establish such results, a basic ingredient will be the weak compactness theory explained as above. This is in fact the main motivation of this note. 

To state our main result, let us introduce some preliminaries about connections on principal bundles; see \cite{Kobayashi-Nomizu-1963-book,Michor-2008-book} for more detailed introduction. First recall that
$Ad(P)$ is the Riemannian vector bundle associated to the principle
bundle $P$ via the adjoint representation $Ad:G\to GL(\mf{g})$.
Denote by $\Ga(M,Ad(P)\otimes\wedge^{k}T^{\ast}M)$ the space of smooth
sections of $Ad(P)$-valued $k$-forms on $M$, and by $\Ga_{c}(M,Ad(P)\otimes\wedge^{k}T^{\ast}M)$
the subspace of compactly supported sections in $\Ga(M,Ad(P)\otimes\wedge^{k}T^{\ast}M)$. The Lie algebra structure
of fibers of $Ad(P)$ induces a natural wedge product $[\al\wedge\be]$
such that, for any $\al\in\Ga(Ad(P)\otimes\wedge^{k}T^{\ast}M)$ and
$\be\in\Ga(Ad(P)\otimes\wedge^{l}T^{\ast}M)$, we have
\[
[\al\wedge\be]|_{U}=[\al_{I},\be_{J}]dx^{I}\wedge dx^{J}
\]
on any chart $(U,x)$ of $M$, where $\al|_{U}=\al_{I}dx^{I}$ and
$\be|_{U}=\be_{J}dx^{J}$ with $I,J$ corresponding multi-indexes.
The inner product structures of $M$ and $\mf{g}$ also induce an
inner product on $Ad(P)\otimes\wedge^{k}T^{\ast}M$, which we still
denote by $\langle,\rangle$ whenever there is no potential confusion.

For any fixed reference connection $\tilde{\na}$, since the space
of connections on $P$ is affine, we define the space of locally weak
$W^{k,p}$-connections on $P$ as
\[
{\cal A}_{\loc}^{k,p}(P)=\left\{ \na_{A}=\tilde{\na}+A:A\in W_{\loc}^{k,p}(M,Ad(P)\otimes T^{\ast}M)\right\} .
\]
Locally on any chart $(U,x)$, if we choose $\tilde{\na}=d$ such
that $\na_{A}|_{U}=d+\ta$ for some $\ta\in W_{\loc}^{k,p}(U,\g\otimes T^{\ast}U)$,
then the curvature form $F_{A}$ on $U$ is represented by
\[
F_{A}|_{U}=d\ta+\frac{1}{2}[\ta\wedge\ta].
\]
Note that $F_{A}\in L_{\loc}^{2}$ whenever $A\in{\cal A}_{\loc}^{0,4}\cap{\cal A}_{\loc}^{1,2}(P)$
(or equivalently $\ta\in L_{\loc}^{4}\cap W_{\loc}^{1,2}$). Each
connection $\na_{A}\in{\cal A}_{\loc}^{k,p}(P)$ induces an exterior
covariant derivative
\[
D_{A}:\Ga(M,Ad(P)\otimes\wedge^{k}T^{\ast}M)\to\Ga(M,Ad(P)\otimes\wedge^{k+1}T^{\ast}M),
\]
such that locally for any section $\psi\in\Ga(U,\g\otimes\wedge^kT^{\ast}U)$,
we have
\[
D_{A}\psi=d\psi+[\ta\wedge\psi].
\]
Accordingly, the adjoint exterior covariant derivative $D_{A}^{\ast}$
with respect to the inner product of $\mf{g}\otimes T^{\ast}M$ can
be explicitly given by
\[
D_{A}^{\ast}=(-1)^{nk+1}\ast D_{A}\ast:\Ga(M,Ad(P)\otimes\wedge^{k+1}T^{\ast}M)\to\Ga(M,Ad(P)\otimes\wedge^{k}T^{\ast}M),
\]
such that locally for any $\psi\in\Ga(U,\g\otimes\wedge^{\ast}T^{k+1}U)$,
there holds
\[
D_{A}^{\ast}\psi=d^{\ast}\psi+(-1)^{kn+1}\ast[\ta\wedge\ast\psi].
\]
Thus we define a weak $\Om$-Yang-Mills connection as follows.

\begin{definition} We say that $A\in{\cal A}_{\loc}^{0,4}(P)$ is
a weak $\Om$-Yang-Mills connection, if $F_{A}$ belongs to $L_{\loc}^{2}(M,Ad(P)\otimes\wedge^{2}T^{\ast}M)$
and
\[
\int_{M}\left\langle F_{A}+\ast(F_{A}\wedge\Om),D_{A}\phi\right\rangle dV_{M}=0
\]
holds for any $\phi\in\Ga_{c}(M,Ad(P)\otimes T^{\ast}M)$. 
\end{definition}

Our main result reads as follows.

\begin{theorem}[Weak compactness]\label{thm: main results} Let
$\tilde{\na}$ be a fixed reference connection on the principle bundle
$P\to M$ and $\na_{A_{i}}=\tilde{\na}+A_{i}\in{\cal A}_{\loc}^{0,4}(P)$
a sequence of weak $\Om$-Yang-Mills connections. Assume that
\[
A_{i}\wto A\qquad\text{in }L_{\loc}^{4}(M,Ad(P)\otimes T^{\ast}M)
\]
and
\[
\sup_{i}\|F_{A_{i}}\|_{L^{2}(E)}<\wq\qquad\text{for any bounded set }E\subset M.
\]
Then $\tilde{\na}+A$ is a weak $\Om$-YM connection as well. 
\end{theorem}

In our follow-up works, we shall apply this weak compactness theorem to extend the Naber-Valtorta theories \cite{Naber-Val-2017-Ann,Naber-Val-2019-Invent} to $\Omega$-YM connections. 

\section{Proof of main result}

In this section we prove Theorem \ref{thm: main results}. First we
prove the following auxiliary lemma.

\begin{lemma}\label{lem: weak continuity of Omega-curvature} Let
$\tilde{\na}$ be a fixed reference connection on the principle bundle
$P\to M$ and $\na_{A_{i}}=\tilde{\na}+A_{i}\in{\cal A}^{0,p}(P)$
is a sequence of weak connections on $P$ for some $p>2$. Suppose
$A_{i}$ satisfies
\[
A_{i}\wto A\qquad\text{in }L_{\loc}^{p}(M,Ad(P)\otimes T^{\ast}M)
\]
for some weak connection $\na_{A}=\tilde{\na}+A\in{\cal A}^{0,p}(P)$,
and
\[
\sup_{i}\|F_{A_{i}}\|_{{\cal M}(E)}<\wq
\]
for all bounded set $E\subset M$, where ${\cal M}(E)$ is the space
of Borel regular real measures on $E$. Then
\[
\lim_{i\to\wq}\int_{M}\langle F_{A_{i}}+\ast(F_{A_{i}}\wedge\Om),\phi\rangle dV_{M}=\int_{M}\langle F_{A}+\ast(F_{A}\wedge\Om),\phi\rangle dV_{M}
\]
for any $\phi\in\Ga_{c}(M,Ad(P)\otimes\wedge^{k}T^{\ast}M)$. \end{lemma}
\begin{proof}
It has been proved by Chen-Giron \cite[Proposition 2.3]{Chen-Giron-21}
that under the above assumptions $F_{A_{i}}\to F_{A}$ in the sense
of distributions. We only need to show the weak continuity of the
second term $\ast(F_{A_{i}}\wedge\Om)$. For any $\phi\in\Ga_{c}(M,Ad(P)\otimes\wedge^{k}T^{\ast}M)$,
we have
\[
\langle\ast(F_{A_{i}}\wedge\Om),\phi\rangle dV_{M}=\phi\wedge\ast\ast(F_{A_{i}}\wedge\Om)=\phi\wedge(F_{A_{i}}\wedge\Om).
\]
Hence
\[
\langle\ast(F_{A_{i}}\wedge\Om),\phi\rangle dV_{M}=F_{A_{i}}\wedge\phi\wedge\Om=\langle F_{A_{i}},\ast(\phi\wedge\Om)\rangle dV.
\]
This gives
\[
\int_{M}\langle\ast(F_{A_{i}}\wedge\Om),\phi\rangle dV_{M}=\int_{M}\langle F_{A_{i}},\ast(\phi\wedge\Om)\rangle dV_{M},
\]
which tends to
\[
\int_{M}\langle F_{A},\ast(\phi\wedge\Om)\rangle dV_{M}=\int_{M}\langle\ast(F_{A}\wedge\Om),\phi\rangle dV_{M}
\]
by the distributional convergence of $F_{A_{i}}\to F_{A}$. This shows
that
\[
\ast(F_{A_{i}}\wedge\Om)\to\ast(F_{A}\wedge\Om)\quad\text{in the sense of distributions.}
\]
The proof is complete.
\end{proof}
The following lemma is a special case of Chen-Giron \cite[Lemma 3.1]{Chen-Giron-21} 
(namely, taking $k=2$ therein), which plays a crucial role in the
proof of Theorem \ref{thm: main results}. 

\begin{lemma}[{Weak continuity, \cite[Lemma 3.1]{Chen-Giron-21}}]\label{lem: Chen-Giron} Let $M_{g}^{n}$
be an $n$-dimensional Riemannian manifold and let $1<p_{1},p_{2}<\wq$,
$1\le\mu_{1},\mu_{2}\le n$, such that
\[
\frac{1}{p_{1}}+\frac{1}{p_{2}}=1,\qquad\mu_{1}+\mu_{2}=s\le n.
\]
Assume that $\{A_{k}^{i}\}_{k\ge1}\subset L_{\loc}^{p_{i}}(M,\mf{g}\otimes\wedge^{\mu_{i}}T^{\ast}M) (i=1,2)$ be $\Om$-YM connections, satisfy
\[
A_{k}^{i}\overset{k}{\longrightarrow}A^{i}\qquad\text{weakly in }L_{\loc}^{p_{i}},
\]
\[
\left\{ dA_{k}^{i}\right\} _{k\ge1}\subset\subset W_{\loc}^{-1,p_{i}}.
\]
Then for any $\phi\in\Ga_{c}(M,\mf{g}\otimes\wedge^{s}T^{\ast}M)$
\[
\lim_{k\to\wq}\int_{M}\langle[A_{k}^{1}\wedge A_{k}^{2}],\phi\rangle dV_{M}=\int_{M}\langle[A^{1}\wedge A^{2}],\phi\rangle dV_{M}.
\]
\end{lemma}

Now we can prove Theorem \ref{thm: main results}. Due to the local
nature of the result, we will assume that all the operations below
are done in a chart. In the general case, it is standard to use a
partition of unit argument, see Chen-Giron \cite{Chen-Giron-21} for
instance.
\begin{proof}
\textbf{Step 1.} Since $A_{i}'s$ are weakly $\Om$-YM connections,
for any $\phi\in\Ga_{c}(M,\mf{g}\otimes T^{\ast}M)$ , we have
\[
\int_{M}\langle F_{A_{i}}+\ast(F_{A_{i}}\wedge\Om),D_{A_{i}}\phi\rangle dV_{M}=0.
\]
Using the formula $D_{A_{i}}\phi=d\phi+[A_{i}\wedge\phi]$, it follows
that
\[
\int_{M}\langle F_{A_{i}}+\ast(F_{A_{i}}\wedge\Om),d\phi\rangle dV_{M}=-\int_{M}\left\langle F_{A_{i}}+\ast(F_{A_{i}}\wedge\Om),[A_{i}\wedge\phi]\right\rangle dV_{M}.
\]
By the triple product identity\footnote{For any $\al^{p}\in\Ga(Ad(P)\otimes\wedge^{p}T^{\ast}M)$, $\be\in\Ga(Ad(P)\otimes\wedge^{q}T^{\ast}M)$
and $\ga\in\Ga(Ad(P)\otimes\wedge^{p+q}T^{\ast}M)$, the following
\textbf{triple product identity} holds:
\[
\boxed{\left\langle [\al^{p}\wedge\be],\ga\right\rangle =(-1)^{p(n-p)}\left\langle \al^{p},\ast[\be\wedge\ast\ga]\right\rangle .}
\]
}, we have
\[
\begin{aligned}\left\langle F_{A_{i}}+\ast(F_{A_{i}}\wedge\Om),[A_{i}\wedge\phi]\right\rangle  & =\left\langle [\phi\wedge A_{i}],F_{A_{i}}+\ast(F_{A_{i}}\wedge\Om)\right\rangle \\
 & =(-1)^{n-1}\left\langle \phi,\ast\left[A_{i}\wedge\ast\left(F_{A_{i}}+\ast(F_{A_{i}}\wedge\Om)\right)\right]\right\rangle \\
 & =\left\langle \phi,\ast\left[\ast\left(F_{A_{i}}+\ast(F_{A_{i}}\wedge\Om)\right)\wedge A_{i}\right]\right\rangle .
\end{aligned}
\]
Hence
\begin{equation}
\int_{M}\langle F_{A_{i}}+\ast(F_{A_{i}}\wedge\Om),d\phi\rangle dV_{M}=-\int_{M}\left\langle \ast\left[\ast\left(F_{A_{i}}+\ast(F_{A_{i}}\wedge\Om)\right)\wedge A_{i}\right],\phi\right\rangle dV_{M}.\label{eq: temp-1}
\end{equation}
Now Lemma \ref{lem: weak continuity of Omega-curvature} implies that
\begin{equation}
\text{LHS of }\eqref{eq: temp-1}\to\int_{M}\langle F_{A}+\ast(F_{A}\wedge\Om),d\phi\rangle dV_{M}.\label{eq: convergence of LHS}
\end{equation}
We need to determine the convergence of the RHS of (\ref{eq: temp-1}).
we will use Lemma \ref{lem: Chen-Giron}. 

\textbf{Step 2.} First note that $A_{i}$ is bounded in $L^{4}$ and
so
\[
\{dA_{i}\}\quad\text{is a bounded sequence in }W^{-1,4}(E).
\]
On the other hand, $dA_{i}=F_{A_{i}}-\frac{1}{2}[A_{i}\wedge A_{i}]$
is bounded in $L^{2}(E)$ by the assumption. This implies that
\[
\{dA_{i}\}\subset\subset W^{-1,2}(E).
\]
As a result, an interpolation argument shows that
\[
\{dA_{i}\}\subset\subset W^{-1,3}(E).
\]

Denote $F_{\Om,i}\equiv F_{A_{i}}+\ast(F_{A_{i}}\wedge\Om)$. Since
$\Om$ is bounded, we have $\ast F_{\Om,i}\in L^{2}(E)$ is a bounded
sequence. Hence
\[
\{d\ast F_{\Om,i}\}\quad\text{is a bounded sequence in }W^{-1,2}(E).
\]
On the other hand, since $A_{i}$ is an $\Om$-Yang-Mills connection,
we have
\[
d^{\ast}F_{\Om,i}=\ast[\ast F_{\Om,i}\wedge A_{i}].
\]
Hence
\[
d\ast F_{\Om,i}=\ast d^{\ast}F_{\Om,i}=(-1)^{n-1}[\ast F_{\Om,i}\wedge A_{i}].
\]
Thus $d\ast F_{\Om,i}$ is bounded in $L^{4/3}(E)$ by H\"older's inequality,
which implies that
\[
\{d\ast F_{\Om,i}\}\subset\subset W^{-1,4/3}(E).
\]
Consequently interpolation arguments imply that
\[
\{d\ast F_{\Om,i}\}\subset\subset W^{-1,3/2}(E).
\]
Now applying Lemma \ref{lem: Chen-Giron} with $p_{1}=3/2$ and $p_{2}=3$,
we derive that
\begin{equation}
\text{RHS of }\eqref{eq: temp-1}\to-\int_{M}\left\langle \ast\left[\ast\left(F_{A}+\ast(F_{A}\wedge\Om)\right)\wedge A\right],\phi\right\rangle dV_{M}.\label{eq: convergence of RHS}
\end{equation}

\textbf{Step 3.} Finally, combining (\ref{eq: convergence of LHS})
and (\ref{eq: convergence of RHS}) yields
\[
\int_{M}\langle F_{A}+\ast(F_{A}\wedge\Om),D_{A}\phi\rangle dV_{M}=0.
\]
 The proof is complete.
\end{proof}

\end{document}